\documentclass[12pt]{amsart}
\usepackage{geometry}          
\usepackage{graphicx}
\usepackage{amssymb}
\usepackage{epstopdf}
\def\lab#1{\label{#1}}

\newtheorem{theo}[subsection]{Theorem}

\newtheorem{lemm}[subsection]{Lemma}
\newtheorem{defi}[subsection]{Definition}
\newtheorem{prop}[subsection]{Proposition}

\newtheorem{rema}[subsection]{Remark}

\def\equat{\refstepcounter{subsection}$$~}

\def\endequat{\leqno{(\boldsymbol{\arabic{section}.\arabic{subsection}})}~$$}

\def\F{{\Bbb{F}}}\def\Q{{\Bbb{Q}}}
\def\C{{\Bbb{C}}}

\def\Z{{\Bbb Z}}

\def\BB{{\bf B}}
\def\CC{{\bf C}}
\def\GG{{\bf G}}\def\GD{{{\bf G}^*}}\def\GF{{{\bf G}^F}}

\def\FF{{\bf F}}
\def\HH{{\bf H}}\def\HF{{\HH^F }}
\def\LL{{\bf L}}\def\LF{{{\bf L}^F}}\def\LD{{{\bf L}^*}}
\def\MM{{\bf M}}
\def\PP{{\bf P}}\def\QQ{{\bf Q}}

\def\SS{{\bf S}}
\def\TT{{\bf T}}\def\TD{{\TT^*}}
\def\UU{{\bf U}}

\def\XX{{\bf X}}
\def\YY{{\bf Y}}

\def\CCC{{\mathcal{C}}}

\def\III{{\mathcal{I}}}

\def\MMM{{\mathcal{M}}}

\def\RRR{{\mathcal{R}}}

\def\VVV{{\mathcal{V}}}

\def\aa{{\alpha}}\def\al{{\alpha}}

\def\dd{{\delta}}

\def\ep{{\epsilon}}
\def\si{{\sigma}}
\def\cc{{\chi}}
\def\th{{\theta}}
\def\la{{\lambda}}

\def\Sy{{\frak S}}

\def\GGGG{{\frak G}}

\def\ooo{{\circ}}
\def\o#1{{{#1}^\ooo}}

\def\ti#1{\widetilde{#1}}

\usepackage{amsmath}
\usepackage{accents}

\def\Aut{{\rm Aut}}
\def\SL{{\rm SL}}\def\GL{{\rm GL}}\def\PGL{{\rm PGL}}\def\SU{{\rm SU}}

\def\Res{{\rm Res}}
\def\Ind{{\rm Ind}}

\def\Ce#1{{{\rm C}_{#1} }}\def\CCe#1{{{\rm C}^\circ_{#1} }}

\def\No#1{{{\rm N}_{#1}}}\def\We#1{{{\rm W}_{#1}}}

\def\opp{{\rm opp}}
\def\ad{{\rm ad}}

\def\Irr{{\rm Irr}}

\def\Prod{\prod}

\def\inn{\subseteq}
\def\nni{\supseteq}

\def\smd{\rtimes}
\def\bull{\hfill\vrule height 1.6ex width .7ex depth -.1ex }
\def\mm{{-1}}

\def\te{\otimes}

\def\ser#1#2{\mathcal{E} (#1 ,#2 )}
\def\Lu#1#2{{\rm R}_{#1}^{#2}}

\def\gen#1{{{<}#1{>}}}
\def\bull{\hfill\vrule height .9ex width .1ex depth -.1ex }
\def\qed{\hfill\vrule height 1.6ex width .7ex depth -.1ex }
\def\mapsot{\ {\leftarrow}{\bull}\ }
\def\lan{\langle}\def\ran{\rangle}
\def\summ{\mathop{\overline{\sum}}}

\def\mr#1{\smash{\mathop{\relbar\joinrel\longrightarrow}\limits^{#1}}}

\def\md#1{\Bigg\downarrow\rlap{$\vcenter{\hbox{$\scriptstyle#1$ }} $}}

\def\Pr#1{Proposition~\ref{#1}}\def\Th#1{Theorem~\ref{#1}}\def\Le#1{Lemma~\ref{#1}}\def\De#1{Definition~\ref{#1}}

\title{On Jordan decomposition of characters for $\SU (n,q)$.}

\author{Marc Cabanes} \address{Institut de Math\'ematiques de Jussieu, Universit\'e Paris 7, 175 rue du Chevaleret,  F-75013 Paris, FRANCE} \email{cabanes@math.jussieu.fr}

\begin{document}

\begin{abstract}
As shown by Bonnaf\'e, a step in proving a Jordan decomposition of characters of finite special linear groups is the parametrization of unipotent characters of centralizers of semi-simple elements in projective linear groups. We show the same kind of result in the case of finite special unitary groups. The proof leads to a mild adaptation of Bonnaf\'e's methods expounded in [B99]. The outcome is a Jordan decomposition of characters compatible with Lusztig's twisted induction.
\end{abstract}

\maketitle

\section*{ Introduction}

Finite groups of Lie type deriving from reductive groups whose center is not connected (e.g. $\SL_n(q)$, $\SU_n(q)$ from $\SL_n(\overline{\F}_q)$) have a representation theory which is less well understood than their connected center counterpart. However they are of great interest for finite group theory since the central coverings of finite simple groups of Lie type are of that kind.

It seems important to establish for those groups a Jordan decomposition of characters with properties similar to the ones known for reductive groups with connected center, for instance commutation with Lusztig's twisted induction of characters (see [L77], [L84]). This was proved for finite special linear groups $\SL_n(q)$ by Bonnaf\'e (see [B99], [B06]~\S~27). As pointed in [B06]~32.1, the missing piece to prove the same statement for finite special unitary groups $\SU_n(q)$ is a parametrization of unipotent characters in (non-connected) centralizers of semi-simple elements in projective unitary groups. This is the main purpose of the present paper (see \Th{JorA}). We essentially review the arguments from the proof of the main statement (7.3.2) of [B99], pointing that many of them apply to a more general situation of a wreath product $\GG =\GG_1\wr A$ whose base group $\GG_1$ is reductive of arbitrary type. 

\medskip

{}

\medskip

{}\medskip

{}

\medskip

{}
\section{Notations and background}

When $(v_i)_{i\in I}$ is a finite family of elements of a vector space over $\Q$, one denotes $$\summ_{i\in I}v_i:=|I|^\mm\sum_{i\in I}v_i .$$

\medskip\noindent{\bf 1.1. Finite groups, central functions.}

When $G$ is a group and $H$ a subgroup (or element) of $G$, one uses the usual notations $\Ce G(H)$, $\No G(H)$ for the centralizer and normalizer, respectively. We define $\We G(H)=\No G(H)/H$ for any subgroup $H\inn G$.

When a group $A$ acts on a set $X$, the subset of $X$ of fixed points is denoted by $X^A$, while for any $x\in X$ its stabilizer in $A$ is denoted by $A_x=\{a\in A\mid a.x=x\}$.

If $G$ is a finite group, one denotes by $\CCC (G)$ the $\C$-vector space of class functions $G\to\C$ endowed with the scalar product $\lan f,f'\ran_G=\summ_{g\in G}f(g)\overline{f'(g)}$. The set of irreducible (ordinary, complex) characters
$\Irr (G)$ is an orthonormal basis of $\CCC (G)$. The action of $\Aut (G)$ on $\CCC (G)$ is denoted by $\si.f(g):=f(\si^\mm (g))$ for any $\si\in\Aut (G)$, $f\in\CCC (G)$, $g\in G$.

If $\phi\colon G \to G$ is a group automorphism, one may form the coset $G\phi$ in the natural semi-direct product $G\smd\gen\phi$. It is $G$-stable and the associated partition of $G\phi$ corresponds with a partition of $G$ into so-called ``$\phi$-classes". One denotes by $\CCC (G\phi)$ the space of $G$-invariant functions $f\colon G\phi\to\C$ endowed with the scalar product $\lan f,f'\ran_{G\phi}= \summ_{g\in G}f(g\phi)\overline{f'(g \phi)}$ for $f,f'\in\CCC (G\phi )$. Note that each element of $\CCC (G\phi )$ is not only $G$-invariant but also $G\smd\gen\phi$-invariant (hence the restriction of an element of $\CCC (G\smd\gen\phi )$) since the image of $g\phi$ by $\phi$ is the same as the image by $g^\mm$ conjugacy.

Denote by ${\bf 1}_{G\phi}\in\CCC (G\smd\gen\phi)$ the characteristic function of the subset $G\phi$.

When $G'\leq G$ is a $\phi$-stable subgroup, one denotes by $\Res^{G\phi}_{G'\phi}\colon\CCC (G\phi)\to\CCC (G'\phi)$ the restriction map and its adjoint by $\Ind^{G\phi}_{G'\phi}\colon\CCC (G'\phi)\to\CCC (G\phi)$ which may be defined by $\Ind^{G\phi}_{G'\phi}(f')(g\phi )=|G|^\mm\sum_xf'(xg\phi (x)^\mm \phi)$ where the sum is over $x\in G$ such that $xg\phi (x)^\mm \in G'$.

\medskip\noindent{\bf 1.2. Wreath products and characters.}

If a group $A$ acts on the group $H$ by group automorphisms, one says that the semi-direct product $H\smd A$ is a wreath product if, and only if, $H$ decomposes as a direct product $H=\prod_{i\in I}H_i$ with $A$ acting on $I$, $a.H_i =H_{a.i}$ for any $a\in A$, $i\in I$, and $a$ induces the identity on $H_i$ as soon as $a.i=i$. Breaking $\prod_{i\in I}H_i$ along $I/A$, those semi-direct products $H\smd A$ are direct products of the situation described below.

Let $X$ a finite set on which acts the finite group $A$, and let $H_0$ be another finite group. Denote by $\MMM (X,H_0)$ the set of maps $X\to H_0$ endowed with its natural group structure. Note that $A$ acts by group automorphisms on $\MMM (X,H_0)$ and that $\MMM (X,H_0)\smd A$ is the typical wreath product usually denoted by $H_0\wr A$. For fixed points we have 

\equat\lab{FixPts} \MMM (X,H_0)^A\cong \MMM (X/A,H_0)\endequat

by the composition with the canonical surjection $X\to X/A$. 

As for conjugacy classes and irreducible characters of $\MMM (X,H_0)$, we have $\Irr (\MMM (X,H_0))\cong\MMM (X, \Irr (H_0))$ by a canonical map. 

So, denoting $H:=\MMM (X,H_0)$, one can identify 

\equat\lab{FixCha} \Irr (H^A)\cong \Irr(H)^A\ \ \ \cc_A\mapsot \cc. \endequat


by $\Irr (\MMM(X,H_0)^A) \cong \Irr(\MMM(X/A,H_0)) \cong\MMM(X/A, \Irr(H_0)) \cong\MMM(X, \Irr(H_0))^A \cong \Irr(\MMM(X,H_0))^A =\Irr (H)^A.$  

For characters in $\Irr (H)^A$, recall the canonical extension to the wreath product $H\smd A$ (see [B99]~2.3.1, [JK81]\S~4.3).

 \begin{defi}\label{CanExt}  {\rm For $\cc\in \Irr(H)$, let 
$$\chi{\smd} A_\cc\in \Irr (H\smd A_\cc)$$ satisfying $(\chi{\smd}A_\cc)(ha)=(\chi{\smd}\gen a)(ha)=\chi _{\gen a}(\pi_a(h))$ for any $h\in H$, $a\in A_\cc$, $\cc_{\gen a}$ defined as in (\ref{FixCha}) and $\pi_a\colon H=\MMM (X,H_0)\to H^{\gen a}=\MMM(X/\gen a,H_0)$ associated to any section of the reduction map $X\to X/\gen a$.

Denote by $\cc{\smd}a$ the restriction of $\cc{\smd}\gen a$ to the coset $Ha\inn H\smd \gen a$.}
\end{defi} 

Let us recall that $\cc{\smd}\gen a$, resp. $\cc{\smd} a$, is denoted as $\ti\cc$, resp $\ti\cc_a$, in [B99] and [B06]~\S 33. 

\begin{rema}\label{CaExRe} {\rm Assume $f\in\CCC (H)^A$ is of the type $\te_{x\in X}f_x$ with $f_x\in\CCC (H_0)$ and $f_{a.x}=f_x$ for each $x\in X$, $a\in A$ (for instance all $f_x$'s equal), so that $f$ is defined by $f((h_x)_{x\in X})=\Prod_{x\in X}f_x(h_x)$. Then, generalizing the above, one may define $f_A\in\CCC (H^A)$ and $f{\smd} A\in\CCC (H\smd A)$ extending $f$ by $f_A((h_x)_{x\in X})=\Prod_{x\in X/A}f_x(h_x)$ and $(f{\smd} A) (ha)=f_{\gen a}(\pi_a (h))$
}\end{rema}

Another interesting property of wreath products $H_0\wr A = H\smd A$ is the simple description of $\CCC (Ha)$ for $a\in A$ (see [B99]~2.2.1). 

\begin{prop}\label{CHa}  {\sl $(\cc{\smd}a)_{\cc\in \Irr (H)^{\gen a}}$ form an orthonormal basis of $\CCC (Ha)$.
 }
\end{prop} 

Clifford theorem (see for instance [N98]~8.9) allows to parametrize $\Irr (H)$ in a very explicit way.

 \begin{defi}\label{IHA}  {\rm Let $M$ be a group on which another group $A$ acts by group automorphisms. One denotes by $\III (M ,A)$ (resp. $\III ^\wedge(M ,A)$) the set of pairs $(\cc ,\al )$ with $\cc\in \Irr (M)$ and $\al\in \Irr (A_\cc )$ (resp. $\al\in A_\cc$). The group $A$ acts on both sets by $\beta .(\cc ,\al )=(^\beta\cc ,^\beta \al )$. 
 
 The corresponding quotients are denoted by $\overline\III (M ,A)$ (resp. $\overline\III ^\wedge(M ,A)$). 
 
 The class of $(\cc ,\al )$ is denoted by $\cc *\al$.}
\end{defi}

 \begin{prop}\label{WrClif}  {\sl In the situation of a wreath product $H\smd A$, one has a bijection $\overline\III (H,A)\mr{\sim} \Irr (H\smd A)$ defined by $\cc *\aa\mapsto \Ind^{H.A}_{H.A_\cc}((\cc{\smd}A_\cc).\al)$.
 }
\end{prop}

\medskip\noindent{\bf 1.3. Non-connected reductive groups}

In this section we give some terminology and basic tools on the representation theory of finite groups of Lie type obtained from non connected reductive groups. The basic reference is [DM94] with generalizations given in [B99]~\S~6.

Note that what follows is used essentially in cases where $\GG =\o\GG\smd A$ where $\o\GG$ (neutral component) is a direct product of general linear groups and $A$ acts by permutations of the summands (``wreath products"). A basic property of this case is that $A$ induces ``quasi-central automorphisms" (see [B99]~5.1.1) of $\GG^\ooo$.

In this section, $(\GG ,F)$ is a (possibly non-connected) reductive group defined over $\F_q$ with associated Frobenius endomorphism $F\colon\GG\to\GG$.

One calls parabolic subgroups of $\GG$ the closed subgroups $\PP\inn \GG$ containing some Borel subgroup of $\o\GG$. Then $\o\PP$ is a parabolic subgroup of $\o\GG$. Denoting by ${\rm R}_{\rm u}$ the unipotent radical, one has $\o\PP ={\rm R}_{\rm u}(\PP)\smd\LL_0$ for $\LL_0$ a so-called Levi subgroup of $\o\GG$ ([DM91]~1.15), and one gets $\PP ={\rm R}_{\rm u}(\PP )\smd \No\PP(\LL_0)$ (see [B99]~\S 6.1).

A group of the above kind $\LL :=\No\PP (\LL_0)$ is called a \it Levi subgroup \rm of $\GG$. A \it quasi-Borel \rm subgroup of $\GG$ is the normalizer in $\GG$ of any Borel subgroup of $\GG$. All quasi-Borel subgroups of $\GG$ are $\o\GG$-conjugate. A \it maximal quasi-torus \rm is any Levi subgroup of a quasi-Borel. It is clearly of the form $\No\GG (\TT_0,\BB_0)$ where $\TT_0\inn\BB_0\inn\o\GG$ are a maximal torus and a Borel subgroup of $\o\GG$. They are all $\o\GG$-conjugate.

In the following, we assume that $a\in\GF$ induces a quasi-central automorphism of $\GG^\ooo$.  There exists a $\gen{a,F}$-stable maximal torus $\TT_0$ of $\o\GG$ ([DM94]~1.36.(ii)).

 \begin{prop}\label{NCTori} \rm ([DM94]~1.40, [B99]~6.1.5). {\sl  Assume $\GG =\o\GG .\gen a$. If $\TT$ is an $F$-stable maximal quasi-torus of $\GG$, it has a $(\o\GG .\gen a )^F$-conjugate $\TT '$ which contains $a$. Moreover $(\TT '{}^{\gen a})^\ooo$ is an $F$-stable maximal torus of $\o{(\GG^{\gen a})}$.
 
 Then $\TT\mapsto  (\TT '{}^{\gen a})^\ooo$ induces a bijection between $(\o\GG .\gen a )^F$-classes of $F$-stable maximal quasi-tori of $\o\GG .\gen a$ and $(\GG^{\gen a})^\ooo{}^F$-classes of $F$-stable maximal tori of $(\GG^{\gen a})^\ooo{}$. In particular the former is in bijection with $F$-classes (see \S~1.1 above) of elements of $\We{\o\GG}(\TT_0)^{\gen a}$. Maximal quasi-tori are called \it of type $w\in\We{\o\GG}(\TT_0)^{\gen a}$ with regard to $\TT_0\gen a$ \sl when they correspond with the $F$-class of $w$ by the above.
 }
\end{prop}

When $\LL$ is an $F$-stable Levi subgroup of $\GG$, the definition of a functor $$\Lu{\LL\inn \PP}\GG\colon \Z\Irr (\LF)\to \Z\Irr (\GF)$$ is similar to the connected case. The variety $\YY =\{ g\in \GG\mid g^\mm F(g)\in {\rm R}_{\rm u}(\PP )\}$ is acted on by $\LF$ on the right and $\GF$ on the left. So the virtual $\overline\Q_\ell$-module obtained from the $\ell$-adic cohomology of $\YY$ with compact support ${\rm H}_{\rm c}(\YY )=\sum_{i\in\Z}(-1)^i {\rm H}_{\rm c}^i(\YY ,\overline\Q_\ell )$ is a virtual $\overline\Q_\ell [\GF\times (\LF{})^\opp ]$-module. Each simple $\overline\Q_\ell [\GF\times (\LF{})^\opp ]$-module can be realized over a finite extension of $\Q$, so we assume the generalized character of $\GF\times (\LF{})^\opp $ we are considering is over $\C$. The associated tensor product functor is the sought $\Lu{\LL\inn\PP}\GG$.

Many properties generalizing the ones of the $\GG =\o\GG$ case hold. 

 \begin{prop}\label{NCRLG}  {\sl Let $\PP$ be a parabolic subgroup of $\GG$ with $F$-stable Levi subgroup $\LL$.
 
 (i) If $\QQ$ is a parabolic subgroup of $\GG$ with $F$-stable Levi subgroup $\MM$ with $\LL\inn\MM$ and $\PP\inn\QQ$, then $$\Lu{\LL\inn \PP}\GG =\Lu{\MM\inn \QQ}\GG \circ \Lu{\LL\inn \PP\cap\MM}\MM .$$
 
 (ii) If $\LL$ is a quasi-maximal torus, the functor $\Lu{\LL\inn\PP}\GG$ does not depend on the choice of $\PP$.
 
 (iii)  If $\cc\in\Irr (\LF)$ and $g\in \GF$, then $$\Lu{\LL\inn \PP}\GG \cc (g)=\sum_x (\Lu{\o\LL \gen{l_x}\inn \o\PP \gen{l_x}}{\o\GG \gen{l_x}}\Res^\LF_{\o\LL^F \gen{l_x}}\cc)\circ\ad_x(g)$$ where the sum is over cosets $x\LL^F.\o\GG^F\inn\GF$ such that $^xg\in\LL^F.\o\GG^F$ and where $l_x\in\LF$ is such that $^xg\in l_x.\o\GG^F$.
 }
\end{prop} 

\noindent{\it Proof.} (i) is [B99]~6.3.3. For (ii) see [DM94]~4.5, [B99]~Remarque~7.1. (iii) is [DM94] 2.3.(i).

\qed

We are mainly concerned with the set $\ser\GF 1$ of so-called ``unipotent" characters (see [DM91]~13.19).

 \begin{defi}\label{NCUnipDef}  {\rm The unipotent characters of $\GF$ are the irreducible components of the various $\Ind^\GF_{\o\GG^F}(\cc )$ for $\cc\in\ser{\o\GG^F}1$. They are also the irreducible components of the generalized characters $\Lu\TT\GG (1)$ for $\TT$ ranging over $F$-stable maximal quasi-tori of $\GG$ (see [B99]~6.4.2).
 }
\end{defi}

\section{ Wreath products of general linear groups}

Adapting [B99]~\S 7, one defines here certain groups $\HH^F$ where $\HH$ is reductive but non-connected in the form $\HH^\circ\smd A$, a wreath product in the sense of \S 2.2 above. We also show that centralizers of semi-simple elements in projective linear or unitary groups are of that type.

Let $q$ be a power of $p$, the characteristic of $\F$.

 \begin{defi}\label{NCtypA}  {\rm
  Let $(\HH ,F)$ be a finite product of pairs $(\HH_i, F_i)$ where $\HH_i=\GL_{n_i}(\F )\wr A_i=\GL_{n_i}(\F )^{d_i}\smd A_i$ ($A_i\inn \Sy_{d_i}$) while $F_i\colon\HH_i\to \HH_i$ is $\si_iF'_i$ for $\si_i\in \No{\Sy_{d_i}}(A_i)$ and $F'_i\colon\HH_i\to \HH_i$ is the map trivial on $A_i$ and raising matrix entries to the $q$-th power composed with some power of the map sending matrices to their transpose-inverse. One denotes by $A$ the product of the $A_i$'s so that $\HH =\HH^\circ\smd A$, and by $\si$ the product of the $\si_i$'s.
}
\end{defi} 

Note that $\HH^F=\o\HH^F\smd A^F$ where $F$ acts by conjugacy by $\si$ on $A$ in $\prod_i\Sy_{d_i}$.

 \begin{rema}\label{QuaCen}  {\rm The action of $A$ on $\o\HH$ is by quasi-central automorphisms (see [B99]~7.1.1).
 }
\end{rema}

One denotes by $\TT_0\inn\o\HH$ the product of subgroups of diagonal matrices in each $\GL_{n_i}(\F )$, i.e. the $F$-stable diagonal torus of $ \o\HH$. One denotes the associated Weyl group $$\o W=\No{\o\HH}(\TT_0)/\TT_0 $$ which identifies with the product of the subgroups of  permutation matrices in each summand $\GL_{n_i}(\F )$. Note that it is stable by $A$ and each element of $\o W$ is fixed by $F'_i$, so that $\si$ is the automorphism of $\o W$ induced by $F$.

 \begin{rema}\label{WrNwr}  {\rm The semi-direct product $\o W\smd A$ is to be considered as the Weyl group of the non-connected reductive group $\HH =\o\HH\smd A$. Note that $(\o W\smd A)^F=\o W^{\gen \si}\smd A^{\gen\si}$ is a wreath product. On the other hand $\HH^F$ may not be one. An example is $\o\HH =(\GL_n(\F ))^d$ with $\si =(1,\dots ,d)$, $F_1$ raising matrix entries to the $q$-th power, $A=\gen\si$. Then $\HH^F=\GL_n(
\F_{q^d})\smd {\gen{F_1}}$. }
\end{rema}

 \begin{prop}\label{WrLev}  {\sl (i) $A$-stable maximal tori of $\o\HH$ are conjugated under $\o\HH{}^A$.
 
 (ii) Suppose $A$ abelian and $\LL$ is an $F$-stable Levi subgroup of $\HH$. Then $(\o\LL .\LF ,F)$ is of the type described in Definition~\ref{NCtypA}.
 }
\end{prop} 

\noindent{\it Proof of~Proposition~\ref{WrLev}.} (i) A maximal torus in a direct product $\o\HH =\prod_{i\in I}\HH^\ooo_i$ of connected groups is the direct product of its intersections with the $\HH_i$'s. A first consequence is that one may assume the $A\inn\Sy_d$ is transitive on a direct product $\o\HH =\HH_1\times\dots \times \HH_1$ ($d$ summands), so that $\o\HH{}^A=\{(g,\dots ,g)\mid g\in\HH_1\}$. An $A$-stable maximal torus will then be of the form $\TT_1\times\dots\times \TT_1$ for $\TT_1$ a maximal torus of $\HH_1$. Then our statement results from conjugacy of maximal tori in $\HH_1$.

(ii) Note first that $\o\LL =\o\HH\cap\LL$, so one can identify $\LL/\o\LL$ with a subgroup $A_\LL$ of $A\cong \HH /\o\HH$. It is $F$-stable. By [B99]~6.2.2, some $\o\HH^F$-conjugate of $\LL$ contains $A_\LL^F$, so one can assume $\o\LL .\LF =\o\LL\smd A_\LL^F$.

Keeping $\HH^\ooo =\prod_{i\in I}\HH^\ooo_i$ and $\o\LL =\prod_{i\in I}\o\LL\cap\HH^\ooo_i$, each $\o\LL\cap\HH^\ooo_i$ is a direct product of $\GL_m$'s. If $\LL\cap\HH^\ooo_{i_0}$ decomposes as $\MM_1\times \dots \times \MM_t$, then for any $a\in A_\LL^F$, $\o\LL \cap \HH^\ooo_{a.i_0}={}^a(\o\LL \cap \HH^\ooo_{i_0})={}^a\MM_1\times\dots\times {}^a\MM_t$ and this decomposition only depends on $a.i_0$, since if $a.i_0=a'.i_0$ then $a^\mm a'$ is trivial on $\HH^\ooo_{i_0}$ and therefore $^a\MM_1 ={}^{a'}\MM_1$, $\dots $ , $^a\MM_t ={}^{a'}\MM_t$. This allows to write $\prod_{i\in A_\LL^F.i_0}\o\LL\cap\HH^\ooo_i$ as a product of $\GL$'s permuted by $A_\LL^F$. Doing it for any $A_\LL^F$-orbit in $I$ gives the same for the whole of $\LL$.

As for the action of $F$, since it is a Frobenius endomorphism, it can be described as a permutation $\si$ of $I$ and a Frobenius twisted or not on each factor $\o\LL_i$ stabilized by a given power $F^j$. The permutations of $I$ induced by $A_\LL^F$ (see above) have to commute with $\si$ since $A_\LL^F$ commutes with $F$ as endomorphisms of $\o\HH$.

This gives our claim.

\qed

\medskip

The relevance of this kind of non-connected groups for the Jordan decomposition of characters is through the following.

Let $n\geq 1$ be an integer. Let $\ti\GG^* =\GL_n(\F)$, let $F^*\colon\ti\GG^*\to\ti\GG^*$ be the raising of matrix entries to the $q$-th power composed with some power of the transpose-inverse automorphism.

Let $\GD =\PGL_n(\F)$ the quotient of $\ti\GG^*$ by its center. \medskip

 \begin{prop}\label{NorLev}  {\sl (i) Let $\LD$ be an $F$-stable Levi subgroup of $\ti\GG^*$, then $(\No{\ti\GG^*} (\LD),F^*)$ is isomorphic with some $(\HH ,F)$ as defined in Definition~\ref{NCtypA}.
 
 (ii) Let $s$ be a semi-simple element in $ (\GG^*)^{F^*}$. Let $\CC^*$ be the inverse image of $\Ce{\GG^*}(s)$ in $\ti\GG^*$. Then $(\CC^*,F^*)$ is isomorphic with some $(\HH ,F)$ as in Definition~\ref{NCtypA} with some cyclic $A\cong \HH /\o\HH$ of order prime to $p$.

 }
\end{prop} 

\noindent{\it Proof of~Proposition~\ref{NorLev}.} (i) We denote by $\TT_0^*$ the torus of diagonal matrices in ${\ti\GG^*}$ so that its normalizer is a semi-direct product of $\TT_0^*$ with the subgroup $\VVV$ of all permutation matrices. Choose as generating set $S\inn\VVV$ the one of transposition matrices exchanging two consecutive elements in the canonical basis and denote by $\LL_I$ the corresponding Levi subgroup for $I\inn S$. The hypothesis is that $\LD ={}^g\LL_I$ for some $g\in {\ti\GG^*}$. Then $(\LD ,F^*)$ is isomorphic with $(\LL_I, xF^*)$ for $x=g^\mm F^*(g)\in \No{\ti\GG^*}(\LL_I)$. But $^x\TT_0^*$ and $\TT_0^*$ are both maximal tori of 
$\LL_I$, so $x=lw$, for $l\in \LL_I$, $w\in \No\VVV (\LL_I)$. Applying Lang theorem to $wF^*$ in $\LL_I$, there is $h\in\LL_I$ such that $h^\mm wF^*(h)w^\mm = l$ and therefore by conjugacy by $h$, $(\LD ,F^*)\cong (\LL_I , xF^*)\cong (\LL_I , wF^*)$. Then $(\No{\ti\GG^*}(\LD ),F^*)\cong (\No{\ti\GG^*}(\LL_I ),wF^*)$. But $\LL_I$ is a product of $\GL_{n_i}(\F )$'s in $\GL_n(\F)$ for $\sum_in_i=n$ and its normalizer is of the type studied before : take $d_{m}$ the number of $i$'s such that $n_i=m$. So that $\No{\ti\GG^*} (\LL_I ) =\LL_I\smd A$ with $A\cong \prod_{m }\Sy_{d_{m}}$ and $w\in A$ with $F^*$ acting trivially on $A$ since $A\inn \VVV$. Hence our claim.

(ii) $\Ce\GD(s)^\ooo$ is an $F^*$-stable Levi subgroup of $\GD$ (see for instance [CE04]~13.15) and $\Ce\GD(s)$ normalizes it, so (i) gives our claim. For the cyclicity and $p'$ order of $\Ce\GD(s)/\Ce\GD(s)^\ooo$, see for instance [CE04]~13.16.(ii).
 
\qed


\section{Parametrizing unipotent characters}

We now take $(\HH ,F)$ as in \De{NCtypA}. Recall $\TT_0$ the torus of diagonal matrices in $\o\HH$ and $\o W\cong \No{\o\HH}(\TT_0)/\TT_0$ the subgroup of permutation matrices. 

As in [B99]~\S 7.4, one wants to parametrize $\ser{\HH^F}1$ (see \S~2.4) by $\Irr (\o W^F.A^F)$. The model is the case of $\o\HH^F$, i.e. finite general linear or unitary groups, as treated by Lusztig-Srinivasan (see [DM91]~\S 15.4). Recall for $w\in\o W$, the torus $(\TT_0)_w$ of $F$-type $w$ in $\o\HH$ (see \Pr{NCTori}). From the case of general linear and unitary groups (Lusztig-Srinivasan [LS77]), corresponding to an irreducible $\o\HH^{\gen\si}$, one easily gets

 \begin{theo}\label{LuSr}  {\rm ([LS77]~2.2, [DM91]~15.8)} {\sl If $\eta\in\Irr (\o W )^F$, let $$R^\ooo_\eta := \summ_{w\in \o W }\eta (w)\Lu{(\TT_0)_w}{\o\HH }(1).$$
 
 Then there is a sign $\ep_\eta\in \{ -1,1\}$ such that $\ep_\eta R^\ooo_\eta\in\Irr (\o\HH^F)$ and one has a bijection 
 
\begin{align*}
\Irr (\o W)^F & \to \ser{\o\HH^F}1  \\
  \eta  &  \mapsto \ep_\eta R^\ooo_\eta
\end{align*}

 }
\end{theo} 

\medskip

\begin{rema}\label{Sign}  {\rm When $\o\HH^F$ is a general linear group, the sign above is 1 (see [DM91]~15.8). When it is a unitary group this sign can be determined from the degree formulas for $R^\ooo_\eta$ (see [Cr85]~13.8) obtained by changing $q\mapsto -q$ in the same formula for general linear groups. It does not depend on $q$. 

The bijection and the formula for $\ep_\eta$ in the general case of $\o\HH^F$ is trivially obtained by direct product. Note however that $\ep_\eta$ takes into account the structure of the summand of $\o\HH^F$, not just of $\o W^F$.}
 
\end{rema} 

\medskip

Looking now for a bijection $\Irr (\o W^F.A^F)\mr{\sim} \ser{\HH^F}1$, one has for the left side a quite simple and canonical parametrization of characters  from the $\o W^F$ case thanks to \Pr{WrClif}. For the right side, one has to find a substitute to the $\cc{\smd} A_\chi$ extension of \De{CanExt} which is not available here since $\HH^F=\o\HH^F\smd A^F$ is in general not a wreath product (see Remark~\ref{WrNwr}). 

Let $\eta\in\Irr (\o W^F)$ and $a\in A^F_\eta$. Having noted that $\o W^F=\o W^\si$ is again a direct product where $a$ acts by permutation of components one has a natural bijection $\Irr (\o W^F)^{\gen a}\to \Irr (\o W^\gen{F,a})\to \Irr (\o W^{\gen a})^F$ (\ref{FixCha}) using the fact that $a$ and $F$ commute. Denote by $\eta_a^F \in \Irr (\o W^{\gen a})^F=\Irr (\o W^{\gen a})^\si$ the image of $\eta$. Since $F$ acts by the permutation $\si$ of components in a decomposition of $\o W^{\gen a}$ as a direct product, one has the associated extension $\eta_a^F{\smd}\gen\si\in\Irr (\o W^{\gen a}\smd\gen\si)$ (see \De{CanExt}).

 \begin{defi}\label{Reta}  {\sl Let $\ti R_\eta$ the element of  $ \CCC (\o\HH^F.A_\eta^F ) $ such that for any $h\in\o\HH^F$, $a\in A_\eta^F$,
 
 $$\ti R_\eta (ha)= \summ_{w\in \o W^{\gen a}}(\eta_a^F{\smd}\gen\si )(w\si )\Lu{(\TT_0.\gen a )_w}{\o\HH .\gen a}(1)(ha),$$ where 
 $(\TT_0.\gen a )_w$ denotes any maximal quasi-torus in $\o\HH\gen a$ of type $w$ with regard to $\TT_0\gen a$ (see \Pr{NCTori}).}

\end{defi} 

\medskip
The fact that $\ti R_\eta$ is indeed a central function on $\o\HH^F.A_\eta^F$ is an easy consequence of the following

 \begin{lemm}\label{bReta}  {\sl  If $b\in A^F$ and $\eta\in\Irr (\o W^F)$, then $^bA_\eta^F=A_{^b\eta}^F$ and $b.\ti R_\eta =\ti R_{^b\eta}$.        }\end{lemm} 

\noindent{\it Proof of~Lemma~\ref{bReta}.}  The first claim is clear. For the second claim, let us fix $a\in A^F_\eta$ and $w\in \o W^{\gen a}$. It suffices to check that $b.\Lu{(\TT_0.\gen a )_w}{\o\HH .\gen a}(1) = \Lu{(\TT_0.\gen {^ba} )_{^bw}}{\o\HH .\gen {^ba}}(1)$. The equivariance of $\ell$-adic cohomology yields $b.\Lu{(\TT_0.\gen a )_w}{\o\HH .\gen a}(1) = \Lu{^b(\TT_0.\gen {a} )_{w}}{\o\HH .\gen {^ba}}(1)$. On the other hand, the fact that $\TT :=(\TT_0.\gen {a} )_{w}$ is of type $w\in\o W^{\gen a}$ in $\o\HH\gen a$ amounts to the fact that there is a $(\o\HH\gen a)^F$-conjugate $\TT '$ of $\TT$ which contains $a$ and such that $(\TT '{}^{\gen a})^\circ$ is of type $w$ in $(\HH^{\gen a})^\circ$, which in turn means that $(\TT '{}^{\gen a})^\circ=((\TT _0{}^{\gen a})^\circ)^x$ for $x\in (\HH^{\gen a})^\circ$ such that $x^\mm F(x)(\TT_0^{\gen a})^\circ = w$ (see \Pr{NCTori}). Then it is clear that $^b\TT$ is  $(\o\HH\gen {^ba})^F$-conjugate to $^b\TT '$ containing $^ba$ and we have $(^b\TT '{}^{\gen {^ba}})^\circ=((\TT _0{}^{\gen {^ba}})^\circ)^{^bx}$ with $(^bx)^\mm F(^bx)(\TT_0^{\gen{^b a}})^\circ = {}^b\bigl( x^\mm F(x)(\TT_0^{\gen a})^\circ\bigr)={}^bw$ since $b$ is $F$-fixed.

 \qed

\medskip

Note that \De{Reta} readily implies $\Res_{\o\HH^F}^{\o\HH^F.A_\eta^F} \ti R_\eta =  R^\ooo_\eta$.

The main result is as follows. The proof is in our last section.

 \begin{theo}\label{NCLuSr}  {\sl Let $\eta\in\Irr (\o W^F)$, then $\ep_\eta\ti R_\eta\in \Irr (\o\HH^F.A_\eta^F)$ 
 
 }
\end{theo} 

\medskip

The consequence on parametrization of unipotent characters of $\HH^F$ is as follows.

 \begin{theo}\label{PaUnNC}  {\sl For $\eta\in\Irr (\o W^F)$, $\xi\in\Irr (A^F_\eta )$, thus defining $ \eta *\xi\in\overline\III (\o W^F ,A^F)$ (see \De{IHA}), let $$R_{\eta *\xi}:= \Ind^{\HH^F}_{\o\HH^F.A_\eta^F}(\xi .\ti R_\eta )\in\CCC (\HH^F).$$ Then one gets a bijection $$ \Irr (\We\HH(\TT_0)^F)
 \to\ser{\HH^F} 1$$ sending $\Ind^{\o W^F\smd A^F}_{{\o W^F\smd A^F_\eta}}(\xi .(\eta{\smd }A_\eta^F))$ to $\ep_\eta R_{\eta *\xi}$.}
\end{theo}

\noindent{\it Proof of~Theorem~\ref{PaUnNC}.} 
We have first $A_\eta^F =A^F_{R^\ooo_\eta}$.
Given the bijectivity of $\eta\mapsto R^\ooo_\eta$ from Theorem~\ref{LuSr}.(ii), this is a consequence of the fact that for any $a\in A^F$, $^aR_\eta^\circ =R^\circ_{^a\eta}$, which is proved in the same fashion as the above \Le{bReta}.

One has $\We\HH (\TT_0)^F\cong \o W^F\smd A^F$, a wreath product, so Clifford theorem in the form of the above \Pr{WrClif} shows that $\eta *\xi\mapsto \Ind^{\o W^F\smd A^F}_{{\o W^F\smd A^F_\eta}}(\xi .(\eta{\smd }A_\eta^F))$ gives a bijection $\overline \III (\o W^F,A^F)\mr{\sim} \Irr (\We\HH (\TT_0)^F)$.

For a given $\ep_\eta R^\ooo_\eta\in\ser{\o\HH^F }1$, $\ep_\eta\ti R_\eta\in\Irr (\o\HH^F.A^F_\eta )$ is an extension by Theorem~\ref{NCLuSr}. We have seen that $\o\HH^F.A^F_\eta $ is the stabilizer of $\ep_\eta R^\ooo_\eta$ in $\HH^F$. Therefore, by Clifford theorem ([N98]~8.9), the irreducible characters of $\HH ^F$ lying over $\ep_\eta R^\ooo_\eta$ are of the type $\ep_\eta R_{\eta *\xi}$ for some $\xi\in\Irr (A_\eta^F)$. This gives a bijection $\eta *\xi\mapsto \ep_\eta R_{\eta *\xi}$ between $\overline \III (\o W^F,A^F)$ and $\ser{\HH^F}1$ by the definition of the latter. 

\qed

\section{Twisted induction}

We keep $(\HH ,F)$, $\TT_0$, $\o W$ as in \S 3. 

In this section, we assume that $A^F$ is abelian and a $p'$-group. This suits the applications we have in mind, thanks to \Pr{NorLev}.(ii).

Here again we follow [B99]~\S 7.6 to give the effect of Lusztig's twisted induction functor on the parametrization of Theorem~\ref{PaUnNC}. The first step is to apply a Mellin transform with regard to the group $A^F$ to our irreducible characters $\ep_\eta R_{\eta *\xi}$. Recall the set $\overline\III^\wedge (\o W^F ,A^F)$ from \De{IHA}.

 \begin{defi}\label{MelTra}  {\rm For $\eta * a\in \overline\III^\wedge (\o W^F ,A^F)$, let $$\hat R_{\eta * a}= \sum_{\xi\in\Irr (A^F_\eta )}\xi (a^\mm )R_{\eta *\xi}.$$
 }
\end{defi} 

Using the commutativity of $A^F$, this is another orthogonal basis of $\C \ser{\HH^F}1$ with base change matrix $({\xi (a)\over |A^F_\cc |}\dd_{\cc ,\cc '})_{\cc *a ,\cc '*\xi}$ and for $h\in \o\HH^F$, $a,b\in A^F$ 

\equat\lab{hatR} \hat R_{\eta *a}(hb)=|A^F_\eta | R_{\eta *1}(hb)\ \ {\rm if} \ \ a=b,\ \ 0\ \ {\rm if} \ \ a\not=b \endequat

see [B99]~7.5.2.

\medskip

Let $\LL$ be an $F$-stable Levi subgroup of $\HH$. By [B99]~6.2.2, one may assume $\LL$ contains $A_\LL^F$ where $A_\LL$ is the image of $\LL$ in $A\cong \HH /\o\HH$. Without changing $\LL^F$, one may even assume that $\LL =\o\LL .\LF =\o\LL .A_\LL^F$, so that $$A_\LL^F =A_\LL =A\cap \LL .$$ Note that this won't change the functor $\Lu{\LL\inn\PP}\HH$ since the parabolic subgroup $\PP$ will be changed into $\o\PP.\LF ={\rm R}_{\rm u}(\PP )\o\LL .\LF$ with same unipotent radical hence same associated variety $\YY$, see \S~1.3. By \Pr{WrLev}.(ii), the parametrization of $\ser\LF 1$ of Theorem~\ref{PaUnNC} applies. 

Let $\TT_\LL$ a diagonal torus of $\o\LL$ which can clearly be taken $A_\LL $-stable. 

 \begin{lemm}\label{ALwF}  {\sl One has $^{g_1} \TT_0=\TT_\LL$ for some ${g_1}\in \o\HH^{A_\LL }$ such that $w_1:={g_1}^\mm F({g_1})\in \o W\cap \o\HH^{A_\LL }$. Moreover $A_\LL^{w_1F}=A_\LL^F =A_\LL =A\cap \LL$. }
\end{lemm} 

\noindent{\it Proof of~\Le{ALwF}.}
Since $\TT_\LL$ is $A_\LL$-stable, \Pr{WrLev}.(i) implies that $^g \TT_0=\TT_\LL$ for some $g\in \o\HH^{A_\LL}$. Then denote $w_1=g^\mm F(g)\in \No\HH (\TT_0)\cap \o\HH^{A_\LL}$. One may also assume that $w_1\in \o W$ since $\No\HH (\TT_0)=\TT_0.\o W$ and if $g^\mm F(g)$ write $t.w_1$, with $t\in \TT_0$, $w_1\in\o W$, one may write $t= t'w_1F(t')^\mm w_1^\mm$ with $t'\in\TT_0$ thanks to Lang theorem applied to $w_1F$ in $\TT_0$, then replace $g$ with $gt'$.

For the last statement, note that, ${g_1}$ being $A_\LL$-fixed and $A_\LL$ being $F$-fixed, $w_1={g_1}^\mm F({g_1})$ is also $A_\LL$-fixed. Then  $A_\LL =A_\LL^F=A_\LL^{w_1F}$.

\qed

Note that $W^\ooo_\LL :=\We{\o\LL}(\TT_\LL)^{g_1}\inn \We{\o\HH}(\TT_0)$ and the pair $(\LL, F)$ is sent to $(\LL^{g_1},w_1F)$ by ${g_1}$-conjugacy. Theorem~\ref{PaUnNC}
for $\LL$ gives a basis $\hat R^\LL_{\la *a}$ (for $\la *a\in\overline{\III}^\wedge (\We{\o\LL} (\TT_\LL)^F, A_\LL )$) of $\C \ser{\LL^F}1$. So determining the effect of the functor $\Lu{\LL\inn\PP}\HH$ on $\ser\LF 1$ amounts to give a formula for $\Lu{\LL\inn\PP}\HH (\hat R^\LL_{\la *a})$. 

Note that the formula does not depend on the choice of $\PP$. Indeed, by \Pr{NCRLG}.(ii), one has independence for tori and, along with transitivity and our definition of $\hat R_{\la *a}$ as a combination of $\Lu\TT\HH$'s, parabolics do not matter. So we now omit $\PP$ in our statements.

Let us fix $\la \in \Irr (\We{\o\LL} (\TT_\LL)^F)$ and $a\in (A\cap \LL)_\la$, thus defining $\la *a\in\overline{\III}^\wedge (\We{\o\LL} (\TT_\LL)^F, A\cap \LL )$. Recall $\la_a^F\in\Irr (\We{\o\LL} (\TT_\LL)^{\gen a})^F$ associated with $\la\in\Irr (\We{\o\LL} (\TT_\LL)^F)^{\gen a}$ thanks to (\ref{FixCha}). Let $\la '_a\in \Irr (W^\ooo_\LL{}^{\gen a})^{w_1F}$ the image of $\la_a^F$ by ${g_1}$-conjugacy $x\mapsto x^{g_1}$ (recall from \Le{ALwF} that ${g_1}$ commutes with $ a$).

Recall the notations from \De{CanExt}.

 \begin{theo}\label{RLGUnNC}  {\sl  One has $$\Lu\LL\HH (\hat R_{\la * a }^\LL)=\sum_{\eta\in\Irr (\o W^F)^{\gen a}}m_\eta \hat R_{\eta * a}^\HH$$ where the scalars $m_\eta\in\C$ are defined by $$\Ind^{\o W^{\gen a}.\si}_{ W ^\ooo_\LL{}^{\gen a}.w_1\si}(\la '_a {\smd}{w_1\si})=\sum_{\eta\in\Irr (\o W^F)^{\gen a}}m_\eta  ( \eta _a{\smd}\si )\ \ \ \  {\rm in}\ \ \CCC (\o W^{\gen a}.\si ).$$}
\end{theo}

Note that the $m_\eta$'s are uniquely defined thanks to \Pr{CHa}.
 
\bigskip

\noindent{\it Proof of~Theorem~\ref{RLGUnNC}.} We essentially recall the main lines of the proof of [B99]~7.6.1. 

First it is clear from (\ref{hatR}) that $\hat R_{\la * a}^\LL$ has support in $\o\LL^F.a$, and therefore also its image by $\Lu\LL\HH$ has support in $\o\HH^F.a$ thanks to the character formula in a case where $|\HH/\o\HH |$ is prime to $p$ (see [DM94]~2.6, 2.7). So it suffices to compare both sides of the claimed equality of Theorem~\ref{RLGUnNC} on restricting to $\o\HH^F.\gen a$.

For the right side, using the definition of $\hat R_{\eta *a}$ and $R_{\eta *\xi}$ (see Theorem~\ref{PaUnNC}), one has $\Res^\HF_{\o\HH^F.\gen a}\hat R^\HH_{\eta *a} = \Res^\HF_{\o\HH^F.\gen a}\Ind^\HF_{\o\HH^F.A^F_\eta}(\sum_{\xi\in\Irr (\o\HH^F.A^F_\eta /\o\HH^F)}\xi (a^\mm )\xi .\ti R_\eta)) $. Moreover, $ \summ_{\xi\in\Irr (\o\HH^F.A^F_\eta /\o\HH^F)}\xi (a^\mm )\xi ={\bf 1}_{\o\HH^F.a}$ (characteristic function of the coset $\o\HH^F.a\inn\o\HH^F.A^F_\eta$). Noting that $\HH^F=\o\HH^F. A^F$ with abelian $A^F$, and using the Mackey formula for the ordinary restriction and induction from subgroups, one gets $\Res^\HF_{\o\HH^F.\gen a}\hat R^\HH_{\eta *a} =\sum_{b\in A^F }\Res^{\o\HH^F.A^F_\eta }_{\o\HH^F.\gen a}{}{}^b({\bf 1}_{\o\HH^F.a}.\ti R_\eta)$. Using now the definition of $\ti R_\eta$ (see \De{Reta}), one gets \equat\lab{RHS} \Res^\HF_{\o\HH^F.\gen a}\hat R^\HH_{\eta *a} = \sum_{b\in A^F}{}^b\bigl({\bf 1}_{\o\HH^F.a}\summ_{w\in\o W^{\gen a}}(\eta_a^F{\smd}\si)(w\si )\Lu{(\TT_0\gen a)_w}{\o\HH\gen a}(1)\bigr)\endequat 

Concerning the left hand side of the claimed equality of Theorem~\ref{RLGUnNC}, applying first \Pr{NCRLG}.(iii), one has $\Res^\HF_{\o\HH^F.\gen a}\Lu\LL\HH (\hat R^\LL_{\la*a})=|A\cap\LL |^\mm \sum_{b\in A^F}{}^b\bigl(\Lu{\o\LL\gen a}{\o\HH.\gen a}(\Res^{\LF}_{\o\LL^F.\gen a} \hat R^{\LL}_{\la*a})\bigr)$. On the other hand, arguments similar to the above allow to write $\Res^{\LF}_{\o\LL^F.\gen a} \hat R^{\LL}_{\la*a} = \sum_{c\in A\cap \LL} {}^c\bigl({\bf 1}_{\o\LL^F.a}\summ_{w\in \We{\o\LL}(\TT_\LL)^{\gen a}}(\la_a^F{\smd}\si )(w\si )\Lu{(\TT_\LL\gen a)_w}{\o\LL\gen a}(1) \bigr)$ and therefore $$\Res^\HF_{\o\HH^F.\gen a}\Lu\LL\HH (\hat R^\LL_{\la*a})=\sum_{b\in A^F}{}^b\bigl(\Lu{\o\LL\gen a}{\o\HH.\gen a}({\bf 1}_{\o\LL^F.a}\summ_{w\in \We{\o\LL}(\TT_\LL)^{\gen a}}(\la_a^F{\smd}\si )(w\si )\Lu{(\TT_\LL\gen a)_w}{\o\LL\gen a}(1))\bigr).$$

By a classical application of the character formula (see [DM91]~12.6 for the connected case) and since $A$ is a $p'$-group, one has $\Lu{\o\LL\gen a}{\o\HH.\gen a}({\bf 1}_{\o\LL^F a}.f)={\bf 1}_{\o\HH^F a}.\Lu{\o\LL\gen a}{\o\HH.\gen a}(f)$ for any $f\in \CCC ({\o\LL^F \gen a})$.
Applying also transitivity of $\Lu\LL\GG$ functors (\Pr{NCRLG}.(i)), one gets 
$\Res^\HF_{\o\HH^F.\gen a}\Lu\LL\HH (\hat R^\LL_{\la*a})=\sum_{b\in A^F}{}^{b}\bigl(  {\bf 1}_{\o\HH^F.a}\summ_{ w\in\We{\o\LL}(\TT_\LL )^{\gen a}}(\la^F_a{\smd} \si)(w\si )\Lu{(\TT_\LL \gen a)_w}{\o\HH\gen a}(1) \bigr)$.

Using now conjugacy by $g_1\in (\o\HH )^{A\cap \LL }$ from \Le{ALwF}, one sees that a maximal quasi-torus of $\o\HH\gen a$ of type $w\in \We{\o\LL}(\TT_\LL )^{\gen a}$ with regard to $\TT_\LL\gen a$ will be of type $w^{g_1}.w_1\in W^\ooo_\LL w_1\inn \o W=\We{\o\HH}(\TT_0)$ with regard to $\TT_0\gen a$. Indeed if that torus is $^g\TT_\LL$ with $g^\mm F(g).\TT_\LL =w$, then $^g\TT_\LL ={}^{gg_1}\TT_0$ with $(gg_1)^\mm F(gg_1)=(g^\mm F(g))^{g_1}.g_1^\mm F(g_1)\in w^{g_1}w_1\TT_0$. Conjugacy by $g_1$ sends $(\We{\o\LL}(\TT_\LL ),a,F)$ to $(W^\ooo_\LL , a, w_1F)$, so one may rewrite now $ \summ_{ w\in\We{\o\LL}(\TT_\LL )^{\gen a}}(\la^F_a{\smd} \si)(w \si )\Lu{(\TT_\LL \gen a)_w}{\o\HH\gen a}(1) =  \summ_{ w'\in W^\ooo_\LL {}^{\gen a}}(\la '_a{\smd} w_1\si)(w'w_1\si )\Lu{(\TT_0 \gen a)_{w'w_1}}{\o\HH\gen a}(1) $, thus implying   

\equat\lab{LHS} \Res^\HF_{\o\HH^F.\gen a}\Lu\LL\HH (\hat R^\LL_{\la*a})= \sum_{b\in A^F}{}^{b}\bigl(  {\bf 1}_{\o\HH^F.a}\summ_{ w'\in W^\ooo_\LL {}^{\gen a}}(\la '_a{\smd} w_1\si)(w'w_1\si )\Lu{(\TT_0 \gen a)_{w'w_1}}{\o\HH\gen a}(1) \bigr)
\endequat

So in view of (\ref{LHS}) and (\ref{RHS}), the claim of our Theorem is established once we can show that 
\equat\lab{L=R} \summ_{ w'\in W^\ooo_\LL {}^{\gen a}}(\la '_a{\smd} w_1\si)(w'w_1\si )\Lu{(\TT_0 \gen a)_{w'w_1}}{\o\HH\gen a}(1) = \endequat

$$\sum_{\eta\in\Irr (\o W^F)^{\gen a}}\! \! m_\eta \! \! \summ_{w\in\o W^{\gen a}}(\eta_a^F{\smd}\si)(w\si )\Lu{(\TT_0\gen a)_w}{\o\HH\gen a}(1) .$$

This follows from a quite formal trick, applying Lemma~3.1.1 in [B99] with (in the notations of that lemma) $H=\o W^{\gen a}$, $K=W^\ooo_\LL{}^{\gen a}$, $x=w_1$, $E(H,\si )=\CCC (\o\HH^F\gen a)$, $\rho^H_w =\Lu{(\TT_0\gen a)_w}{\o\HH\gen a}(1)$ for $w\in H$, $E(K,x\si )=\CCC (\o\LL^F\gen a)$, $\rho^K_{w'} =\Lu{(\TT_0 \gen a)_{w'w_1}}{\o\LL\gen a}(1)$ for $w'\in K$, and $R^H_K =\Lu{\o\LL\gen a}{\o\HH\gen a}$.

\qed

\bigskip

Let $\GG :=\SL_n(\F )$ and $F\colon\GG\to\GG$ the raising of matrix entries to the $q$-th power composed with some power of the transpose-inverse automorphism. Let $\GD =\PGL_n(\F)$ endowed with the Frobenius endomorphism $F^*$ defined in the same way. If $s\in\GD^{F^*}$ is a semi-simple element, recall the rational series $\ser\GF{[s]}\inn\Irr (\GF)$ ([B06]~11.A).
According to [B06]~32.1, the bijection $\ser\GF{[s]}\to\Irr (W(s)^{w_sF^*})$ ([B06]~23.7) and the above \Th{PaUnNC} and \Th{RLGUnNC} (generalizing [B99]~7.4.3 and 7.6.1, respectively) imply the following

 \begin{theo}\label{JorA}  {\sl Assume $q$ is such that the conjecture $\GGGG$ of [B06] \S~14.E is satisfied. 
 
 For any semi-simple element $s\in\GD^{F^*}$, one has a bijection $\aleph_{\GG ,s}\colon \ser\GF{[s]}\to\ser{\Ce\GD(s)^{F^*}}1$, such that for any Levi $\LD$ such that $\LD^{F^*}$ contains $s$ we have a
commutative square }
 \[ \begin{array}{ccc}
 \Z\ser\LF{[s]} &  \mr{\aleph_{\LL ,s}} &  \Z \ser{\Ce\LD(s)^{F^*}}1\\
 {\ep_\GG\ep_\LL\Lu\LL\GG }\md{}\ \ \ \ \ \  \ \ \ \ &   & \md{\ep_{\Ce\GD (s)^\ooo}\ep_{\Ce\LD (s)^\ooo} \Lu{\Ce\LD (s)}{\Ce\GD (s)} }\\
  \Z\ser\GF{[s]} &  \mr{\aleph_{\GG ,s}} &   \Z \ser{\Ce\GD(s)^{F^*}}1
\end{array} \]  

\medskip

\noindent (where we have extended linearly the bijections $\aleph$).
\end{theo}

\noindent{\it Proof.} Thanks to \Pr{NorLev}.(ii), $\HH :=\Ce\GD (s)$ is of the type studied in the preceding section. 

Let's recall why $\Ce\LD (s)$ is a Levi subgroup of $\HH$ in the sense of our \S~1.3. Indeed, if $\PP^* =\UU^*.\LD$ is a parabolic subgroup of $\GD$ of which $\LD$ is a Levi subgroup, the classical description of $\Ce\GD (s)$ in terms of roots (see [DM91]~2.3) implies that $\CCe{\PP^*}(s)$ is a parabolic subgroup of $\Ce\GD (s)^\ooo$ with Levi decomposition $\Ce{\UU^*}(s)\CCe{\LL^*}(s)^\ooo$. Then $\Ce{\PP^*}(s)$ is parabolic subgroup of $\Ce\GD (s)$ and $\Ce\LD (s)=\No{\Ce{\PP^*}(s)}(\Ce{\PP^*}(s)^\ooo,\Ce{\LL^*}(s)^\ooo)$ since the intersection of the latter with $\Ce{\UU^*}(s)$ is trivial, due to $\Ce{\UU^*}(s)=\Ce{\UU^*}(s)^\ooo$ (see [DM91]~2.5) and the Levi decomposition of $\Ce{\PP^*}(s)^\ooo$.

Note that the vertical arrows do not mention the parabolic subgroups. For the one on the right, this is due to the fact mentioned above before stating \Th{RLGUnNC}. Then one may choose matching parabolics on both sides and,    once the theorem is proved, indepedence on the right side implies independence on the left. However, independence of $\Lu{\LL\inn\PP}\GG$ with regard to $\PP$ is known for all groups of classical type as a consequence of the main theorem in [BM11].\def\FD{{F^*}}

Let $A_\GD (s):=\Ce\GD (s)/\CCe\GD (s)$, so that $\Ce\GD (s)=\CCe\GD (s)\smd A_\GD (s)$, both factors being $\FD$-stable (see \Pr{NorLev}.(ii)) and denote $W(s):=\We{\Ce\GD (s)}(\TD )\nni W^\ooo (s):=\We{\CCe\GD (s)}(\TD)$ where $\TD$ is any maximally split torus of $\CCe\GD (s)$. As in [B06]~\S 23, let $w_s\in W^\ooo(s)^{A_\GD (s)^\FD}$ such that $w_s\FD$ acts on $W^\ooo (s)$ by permutation of irreducible components (see [B06]~23.1). 

One defines the bijection $\aleph_{\GG ,s}$ as a composite of two bijections $$\aleph_{\GG ,s}\colon \ser\GF{[s]}\leftarrow\Irr (W(s)^{w_s\FD})\to\ser{\Ce\GD(s)^{F^*}}1.$$ For the left side, one has the isometry $R[s]\colon \C\Irr (W(s)^{w_s\FD})\to \C\ser\GF{[s]}$ of [B06]~23.5 such that for any $\th\in\Irr (W(s)^{w_s\FD})$,  $R[s](\th )\in\ep_\th\ser\GF{[s]}$ for a sign $\ep_\th =\ep_\GG\ep_{\CCe\GD (s)}\ep_\eta$ if $\th$ lies above some $\eta\in\Irr (W^\ooo (s)^{w_s\FD})$ (see [B06]~23.9, this is where the conjecture $\GGGG$ is assumed). The sign $\ep_\eta$ is the one denoted that way in the above \Th{LuSr} as can be seen by applying [B06]~23.13 and 23.15 with $s=1$.

For the right hand side, note first that $W(s)^{w_s\FD}\cong  \We{\Ce\GD (s)}(\TT_0^* ) ^\FD $ where $\TT_0^*$ is a diagonal torus of $\CCe\GD (s)$. This is because, by [B06]~\S 23, the type of $\TT_0^*$ with regard to $\TD$ is the longest element of the Weyl group on each summand of unitary type (and trivial otherwise), which coincides with $w_s$ except for summands of rank 1, but then the fixed point groups are the same. Now the bijection of \Th{PaUnNC} for $\HH :=\Ce\GD (s)$ gives a bijection $\Irr (W(s)^{w_s\FD})\to\ser{\Ce\GD(s)^{F^*}}1$. 

The composite $\aleph_{\GG ,s}$ we get is by $$\ep_\GG\ep_{\CCe\GD (s)}\ep_\eta R[s](\th)\mapsot \th :=\eta *\xi\mapsto\ep_\eta R_\th$$ (notations of \Th{PaUnNC}). 

In view of the formula of compatibility with $\Lu\LL\GG$ functors to prove, we can remove the signs and study the maps extended by linearity $\C\ser\GF{[s]}\leftarrow \C\Irr (W(s)^{w_s\FD})\to\C\ser{\Ce\GD(s)^{F^*}}1$ induced by $ R[s](\th)\mapsot \th :=\eta *\xi\mapsto R_\th$. Denote $A:= A_\GD (s)^\FD$. Each map is a direct sum of the maps $\C\ser\GF{[s]}\leftarrow\CCC (W^\circ (s)^{w_s\FD}a)^A\to\C\ser{\Ce\GD (s)^\FD}1$ defined for each $a\in A_\GD (s)^\FD$ as follows. On the right this is $\eta{\smd}a\mapsto \hat R_{\eta *a}$ from \De{MelTra}. The compatibility we expect is by our \Th{RLGUnNC} for this side and each $a\in A_\LD (s)^F$.

On the left side the map is $R[s,a]\colon \CCC (W^\circ (s)^{w_s\FD}a)^A\to\C\ser\GF{[s]}$ from [B06]~\S 23.C. Through the identification $\CCC (W^\circ (s)^{w_s\FD}a)\cong \CCC (W^\circ(s)^{\gen a}w_s\FD) =\CCC (W^\circ(s)^{\gen a}\FD)$ from the wreath product structure, this is deduced from a map $\RRR [s,a]\colon \CCC (W^\circ(s)^{\gen a}\FD)^A\to\C\ser\GF{[s]}$ (see [B06]~17.18). The compatibility with $\Lu\LL\GG$ is [B06]~17.24. To sum up, the $\RRR[s,a]$ maps for $\GG$ and $\LL$ transform the $\Lu\LL\GG$ functor into an induction of central functions on right cosets in $W(s)$ matching the one satisfied by $\hat R$ seen above. This gives our claim

\qed

\bigskip

{}

\section{Proof of Theorem~\ref{NCLuSr}}

We now give the proof of Theorem~\ref{NCLuSr}. It follows the same steps as in [B99]~\S 8, and we summarize the main ideas. In order to conform to the notations in [B99], we switch back to denoting by $\GG$ (instead of $\HH$) a wreath product of general linear groups as in Definition~\ref{NCtypA}. We keep the other notations $n_i$, $d_i$, $A_i$, $\TT_0$, $\o W$, etc..

Let us recall the construction.

Let $\eta\in\Irr(\o W^F)$, and $\eta_a^F\in\Irr (\o W^{\gen a})$ for $a\in A_\eta$ associated with $\eta$ through the identification $\Irr (\o W^F)^{\gen a} = \Irr (\o W^{\gen{a,F}} )= \Irr (\o W^{\gen a})^F$ (see (\ref{FixCha}). Let $\ti R_\eta^\GG\in\CCC (\o\GG^F.A_\eta^F)$ be defined for $g\in\o\GG^F$, $a\in A^F_\eta$, by (see \De{Reta})

 $$\ti R_\eta^\GG (ga)= \summ_{w\in \o W^{\gen a}}(\eta_a^F{\smd}\gen\si )(w\si )\Lu{(\TT_0.\gen a )_w}{\o\GG .\gen a}(1)(ga),$$ where 
 $(\TT_0.\gen a )_w$ denotes any maximal quasi-torus in $\o\GG\gen a$ of type $w$ with regard to $\TT_0\gen a$ (see \Pr{NCTori}).

 One must prove $\ti R_\eta^\GG\in \pm\Irr (\o\GG^F.A_\eta^F)$.

\bf Step 1. \sl   $\ti R^\GG_1=1$.  \rm 

This is [DM94]~Proposition~4.12.

\bf Step 2. \sl   $\ti R_\eta$ has norm 1.  \rm 

The proof of [B99]~8.1.2, via computation of  scalar products \newline $\lan\Lu{(\TT_0.\gen a )_v}{\o\GG .\gen a}(1) , \Lu{(\TT_0.\gen a )_w}{\o\GG .\gen a}(1)\ran_{\o\GG^F.a}$ (see [DM94]~4.8), applies without change.\qed

\medskip

It now suffices to prove that $$\ti R_\eta\in\Z\Irr (\o\GG^F.A_\eta^F )$$ to have our claim.

By direct product, one assumes that there is just one $i$, so that $\GG^\ooo =\GL_n(\F )^d$. In view of the claim and the definition of $\ti R_\eta$, one may also assume that $A=A^F$, i.e. that $A\inn\Sy_d$ centralizes $\si$. One may also assume $A=\Ce{\Sy_{d}}(\si )$ since this case will imply the same for any subgroup. 
Again by direct product, one assumes that all cycles of $\si$ have same length a divisor $e$ of $d=ek$. 

Then $$\GG =\GL_n(\F )^{ek}\smd A \ \ {\rm with}\ \ A\cong (\Z/e\Z)^{k}\smd \Sy_{k}$$ with $A$ embedding in $\Sy_{ek}$ as $\Ce{\Sy_{ek}}(\si )$ for $\si =(1,\dots , e)(e+1,\dots , 2e)\dots (ek-e+1,\dots ,ek)$.

\bf Step 3. \sl   Assume moreover $k=1$, that is $e=d$. Then $\ti R_\eta$ is a generalized character.  \rm 

Note that $A$ fixes any element hence any character of any subgroup of $\o W^F$.

Let $\LL_1\inn \GL_n(\F)$ be a product of linear groups associated with a composition of $n$ (that is a standard Levi subgroup for the usual BN-pair). Then $\LL:= (\LL_1)^e\smd A$ is an $F$-stable Levi subgroup of $\GG$. Let $W^\ooo_\LL =\o W\cap\LL$ and $\la\in\Irr (W_\LL^\ooo{}^F)$. This is consistent with the notation of \Th{RLGUnNC} above. Then the above process defines a central function $\ti R_{\la}^\LL$ on $\LL^F$.

\begin{lemm}\label{RIndW}  {\sl $$\Lu\LL\GG (\ti R_{\la}^\LL)=\sum_{\eta \in\Irr (\o W^F)}m_\eta  \ti R_{\eta }^\GG$$
 where $m_\eta = \lan \eta  ,\Ind_{W_\LL^\ooo{}^F}^{\o W^F}(\la)\ran_{\o W^F}$ \ }
\end{lemm} 

Let us show why this Lemma is enough to establish (E) in this case. 

By a classical theorem about characters of symmetric groups (see for instance [JK81]~2.2.10), the integral matrix $(\lan \Ind_{W_\LL^\ooo{}^F}^{\o W^F}(1_{W_\LL^\ooo{}^F} ),\eta \ran_{\o W^F})_{\LL ,\eta }$ for $\LL_1$ ranging over Levi subgroups defined by partitions of $n$ and $\eta \in\Irr (\o W^F)$ with $\o W^F\cong \Sy_n$ is invertible. Applying the Lemma~\ref{RIndW} with each $\la$ the trivial character of ${ W^\ooo_\LL{}^F}$, allows to express $\ti R_\eta$ as an integral combination of central functions of type $\Lu\LL\GG (\ti R_1^\LL)$. But Step~1 tells us that $\ti R_1^\LL =1_{\LL^F}$. Hence our claim since $\Lu\LL\GG$ maps generalized characters to generalized characters. 

\noindent{\it Proof of~Lemma~\ref{RIndW}.} One checks the claimed equality on $ga$ with $a=\si^i$ for $i\geq 0$.

The situation is a special case of the one described by \Th{RLGUnNC} with $A=A_\LL$, $w_1=1$ and one has basically to prove the same statement as (\ref{L=R}) with our $m_\eta$'s. Applying [B99]~3.1.1 again, one gets $$\Lu\LL\GG(\ti R_\la^\LL)=\sum_{\eta\in\Irr (\o W^F)}\lan\eta_{\si^i}{\smd}\si ,\Ind^{\o W^{\gen {\si^i}}.\si}_{W_\LL^\ooo{}^{\gen {\si^i}}.\si}(\la_{\si^i}{\smd}\si )\ran_{\o W^{\gen {\si^i}}.\si}\ti R_\eta^\GG .$$ 

There remains to check that $\lan\eta_{\si^i}{\smd}\si ,\Ind^{\o W^{\gen {\si^i}}.\si}_{W_\LL^\ooo{}^{\gen {\si^i}}.\si}\la_{\si^i}{\smd}\si\ran_{\o W^{\gen {\si^i}}.\si}=$ $\lan \eta  ,\Ind_{W_\LL^\ooo{}^F}^{\o W^F}(\la)\ran_{\o W^F}$. This is a general property relating induction on right cosets and ordinary induction in wreath products, see [B99]~3.2.1.

\qed

 \medskip

\bf Step 4. \sl   Assume arbitrary $k$.   Then $\ti R_\eta$ is a generalized character.  \rm 
 
 \medskip
 
 We define a more general framework to which [B99]~\S~8.4 applies.
 
 Take $(\GG_0,F_0)$ a connected reductive group defined over $\F_q$ having an $F_0$-stable maximal torus $\SS_0$ such that $F_0$ acts trivially on $\We{\GG_0}(\SS_0)$. 

Let $$\GG =\GG_0^{ek}\smd A \ \ {\rm with}\ \ A=(\Z/e\Z)^{k}\smd C$$ for $C\inn \Sy_{k}$ embedded as before in $\Sy_{ek}$ as a subgroup of the centralizer $\cong\Z /e\Z\wr \Sy_{k}$ of $\si =(1,\dots , e)(e+1,\dots , 2e)\dots (e(k-1)+1,\dots ,ek)$. Let $F_0$ extend to $\GG$ by acting the same on each factor of $\GG^\ooo = \GG_0^d$ and trivially on $A$, let $F=\si F_0$.

Denote $B\cong (\Z/e\Z)^k$ so that $A\cong B\smd C$. Note that $B$ fixes any element hence any character of $\o W^F$.

Denote $\HH :=\o\GG\smd B =(\GG_0^e\smd \Z/e\Z)^{k}$, so that $\GG =\HH\smd C$ with $F$ acting trivially on the second term and the pair $(\HH ,F)$ being a product of $k$ identical pairs $(\HH_1,F)$ as in the above Step~3. 

Denote $\o W=\We{\o\GG}(\TT_0)$ where $\TT_0=(\SS_0)^d$.

For $\eta\in\Irr (\o W^F)^C=\Irr (\o W^F)^A$, let $\ti R_\eta^\GG \in\CCC (\GF)$ defined by $\Res^\GF_{\o\GG^F.a}\ti R_\eta^\GG =\summ_{w\in\o W^{\gen a}}(\eta_a^F{\smd}a)(w\si) .\Res_{\o\GG^F.a}^{\o\GG^F.\gen a}\Lu{(\TT_0.\gen a)_w}{\o\GG .\gen a}(1)$ for any $a\in A=A_\eta^F$, and $\ti R_\eta^\HH\in\CCC (\HF)$ by $\Res^\HF_{\o\GG^F.b}\ti R_\eta^\HH =\summ_{w\in\o W^{\gen b}}(\eta_b^F{\smd}b)(w\si) .\Res_{\o\GG^F.b}^{\o\GG^F.\gen b}\Lu{(\TT_0.\gen b)_w}{\o\GG .\gen b}(1)$ for any $b\in B=B_\eta^F$. 

Note that the group $\GF$ is also a wreath product $\HF\smd C ={(\HH_1)}^F\wr C$. With this description, $\ti R^\HH_\eta$ is a tensor product of similarly defined representations $\ti R_{\eta_j}^{\HH_j}$ ($j=1,\dots ,k$) of summands $\cong (\HH_1)^F$ with $\eta_j=\eta_{c.j}$ for each $c\in C$. Remark~\ref{CaExRe} then allows to define $\ti R_\eta^\HH{\smd}C\in\CCC (\HF\smd C)$.  

\medskip

 \begin{theo}\label{RetaGH}  {\sl On $\GF = \HH^F\smd C$, one has $\ti R_\eta^\GG = \ti R_{\eta}^{\HH}{\smd}C $.
 }
\end{theo} 

Let's say how this would complete our last step. 

With $\GG_0=\GL_n(\FF)$, $F_0$ the raising of matrix entries to the $q$-th power possibly composed with the transpose-inverse automorphism, $\SS_0$ the diagonal torus, and $C=(\Sy_k)_\eta$ we cover our initial situation. \Th{RetaGH} then implies that $\ep_\eta\ti R_\eta^\GG = \ep_\eta\ti R_{\eta}^{\HH}{\smd}C $ which in turn is an irreducible character thanks to \De{CanExt} and the fact that $\ep_\eta\ti R_{\eta}^{\HH}$ is an irreducible character by Step~3 as said before.

\medskip

\noindent{\it Proof of~Theorem~\ref{RetaGH}.}  
By direct product, one may assume that $C$ is transitive. Given the definition of $\ti R_\eta^\GG$ and the one in Remark~\ref{CaExRe}, one may also assume that $C$ is cyclic. Since it commutes with $\si = (1,\dots , e)(e+1,\dots , 2e)\dots (e(k-1)+1,\dots ,ek)$ in $\Sy_{ek}$, one may then index 
elements in $\{1,\dots , e\},\{e+1,\dots , 2e\}\dots \{e(k-1)+1,\dots ,ek\}$ so that $C$ is generated by $c=(1,\dots ,ke)^e=(1,1+e,\dots ,1+e(k-1))(2,2+e,\dots ,2+e(k-1))\dots  $ and one evaluates the equality on $\HF.c$.

Up to some adequate conjugacy inside $(\Z/e\Z)\wr C$ (see [B99] p. 95), not involving the type of $\GG_0$, one may content ourselves with evaluating our functions on $h^\circ .bc$ where $h^\circ =(h^\circ_1,\dots ,h^\circ_k)\in \o\GG^F=((\GG^e_0{})^F)^k\cong (\GG^{F^e}_0)^k$, $b=(\tau ,1,\dots ,1)^m\in B$ for $\tau$ corresponding with the cycle $(1,\dots ,e)$ and $m\geq 1$.

On the left side of our claim we have \equat\lab{tRG} \ti R_\eta^\GG (\o h .bc)=\summ_{w\in\o W^\gen{bc}}\eta_{bc}^F{\smd}\si(w\si )\Lu{(\TT_0\gen{bc})_w}{\o\GG\gen{bc}}(1)(\o h.bc). \endequat

For the right side one must compute $(\ti R_\eta^\HH{\smd} C)((\o h.b)c)=(\ti R_\eta^\HH){}_c(\pi_c(\o h.b))$ (see \De{CanExt}). 

For $x\in \HH_1=\GG_0^e\smd\Z/ e\Z$, denote $x^{(k)}=(x,\dots ,x)$ ($k$ times), as element of $\HH$. We have $(\HF)^{\gen c}=\{ g^{(k)}  \mid g\in \HH_1^F\}$ and $\pi_c(h_1,h_2,\dots ,h_k) = (h_k\cdots h_1)^{(k)} $ for $h_1,h_2,\dots ,h_k\in \HH_1^F$. So $\pi_c(\o h.b)=( h^\ooo_k\dots  h^\ooo_1\tau^m)^{(k)}=( h^\ooo_k\dots  h^\ooo_1)^{(k)}\si^m$. We have $\ti R_\eta^\HH = (\ti R^{\HH_1}_{\eta_1})^{\te k}$, so $(\ti R_\eta^\HH)_c$ may be identified with $\ti R_{\eta_1}^{\HH_1}$ or $\ti R_{\eta^{\gen c}}^{\HH^{\gen c}}$ where $\eta^{\gen c}(w^{(k)})=\eta_1(w)$ for any $w\in \o W\cap\HH_1^F$. So $(\ti R_\eta^\HH{\smd}C)((h^\ooo   b)c)=\ti R^{\HH^{\gen c}}_{\eta^{\gen c}}(( h^\ooo_k\dots  h^\ooo_1\tau^m)^{(k)}) =\summ_{w\in\o W^{\gen{c,b}}}((\eta^{\gen c})^F_{\si^m}{\smd }\si )(w\si)\Lu{\TT '_w}{(\o\GG\gen{\si^m})^{\gen c}}(1)(( h^\ooo_k\dots  h^\ooo_1)^{(k)}\si^m)$ where $\TT '_w$ denotes the maximal quasi-torus of $(\o\GG\gen{\si^m})^{\gen c}$ of type $w\in \o W^{\gen{c,b}}$.

Comparing this with the formula (\ref{tRG}) above, one notices that $\o W^\gen{bc} =\o W^{\gen{c,b}} $ (since $\gen{bc}$ and ${\gen{c,b}}$ have same orbits on $\{1,\dots ,ek\}$), and similarly $\eta_{bc}^F{\smd}\si = (\eta^{\gen c})^F_{\si^m}{\smd }\si $, which leads to ask if 

\equat\lab{RTGfix} \Lu{(\TT_0\gen{bc})_w}{\o\GG\gen{bc}}(1)(\o h.bc) = \Lu{\TT '_w}{(\o\GG\gen{\si^m})^{\gen c}}(1)(\pi_c(\o h.b)). \endequat

As in [B99]~8.4.3, one first notes that $({(\TT_0\gen{bc})_w})^{\gen c}=\TT '_w$ from equalities in $B\smd C$ and then that (\ref{RTGfix}) holds thanks to 
[DM94]~5.3 relating $\Lu\TT\GG(1)$'s with varieties of type $\XX_w$ (see for instance [CE04]~7.13) and [B99]~8.4.5 giving an isomorphism between varieties $\XX_w$ for $\GG^{\gen c}$ and some corresponding subvarieties for $\GG$, both holding for a $\GG_0$ of any type.
This completes the proof of \Th{RetaGH}.

\qed




\medskip {}


\bigskip

\section*{References.}

\bigskip{}


[B99] C. Bonnaf\'e, Produits en couronne de groupes lin\'eaires, \sl J. Algebra \bf 211 \rm (1999), 57--98.

[B06] C. Bonnaf\'e, Sur les caract\`eres des groupes r\'eductifs finis \`a centre non connexe~: applications aux groupes sp\'eciaux lin\'eaires et unitaires, \sl Ast\'erisque \bf 306\rm , 2006.

[BM11] C. Bonnaf\'e and J. Michel, Computational proof of the Mackey formula for $q > 2$, \sl J. Algebra \bf 327 \rm (2011), 506--526.

[CE04] M. Cabanes and M. Enguehard, {\it Representation theory of finite reductive groups}, Cambridge, 2004.

[Cr85] R. Carter, {\it Finite groups of Lie type : conjugacy classes and complex characters}, Wiley, New York, 1985.

[DM91] F. Digne and J. Michel, {\it Representations of finite groups of
Lie type}, Cambridge, 1991.

[DM94] F. Digne and J. Michel, {Groupes r\'eductifs non connexes} \sl Ann. Sci. \'Ec. Norm. Sup., \rm 4\`eme s\'erie \bf 27 \rm (1994), 345--406.


[JK81] G. James and A. Kerber, {\it The representation theory of the symmetric group}, Addison-Wesley, Reading, 1981.

[L77] G. Lusztig, Irreducible representations of finite classical groups, \sl Inventiones Math. \bf 43 \rm (1977), 125--175.

[L84] G. Lusztig, Characters of reductive groups over a finite field, \sl Ann. of Math. Studies \bf 107\rm , Princeton, 1984.


[LS77] G. Lusztig, B. Srinivasan, The characters of the finite unitary groups, \sl J. Algebra, \bf 49 \rm (1977), 167-171.




[N98] G. Navarro, {\it Characters and blocks of finite groups,} Cambridge, 1998.

\end{document}